\documentclass{ifacconf}

\usepackage{graphicx}      % include this line if your document contains figures
\usepackage{natbib}        % required for bibliography

\usepackage{amsmath}
\usepackage{mathtools}
\usepackage{amssymb}
\usepackage{enumitem}
\newcommand{\ie}{i.\,e.~}
\newcommand{\eg}{e.\,g.~}
\newcommand{\Span}[1]{\mathrm{span}\{#1\}}
\newcommand{\Rank}[1]{\mathrm{rank}(#1)}
\newcommand{\Dim}[1]{\mathrm{dim}(#1)}
\newcommand{\C}[1]{\mathcal{C}(#1)}
\newcommand{\D}{\mathrm{d}}
\newcommand{\Lie}{\mathrm{L}}

\newcommand{\Mod}{~\mathrm{mod}~}

\begin{document}
	\begin{frontmatter}
	
	\title{On a Flat Triangular Form Based on the Extended Chained Form} 
%	% Title, preferably not more than 10 words.
%	
	\thanks[footnoteinfo]{The first author and the second author have been supported by the Austrian Science Fund (FWF) under grant number P 32151 and P 29964.}
	
	\author[First]{Conrad Gst{\"o}ttner} 
	\author[First]{Bernd Kolar} 
	\author[First]{Markus Sch{\"o}berl}
	
	\address[First]{Institute of Automatic Control and Control Systems Technology, Johannes Kepler University, Linz, Austria.\\(e-mail: \{conrad.gstoettner,bernd.kolar,markus.schoeberl\}@jku.at)}

	\begin{abstract}                % Abstract of not more than 250 words.
		In this paper, we present a structurally flat triangular form which is based on the extended chained form. We provide necessary and sufficient conditions for an affine input system with two inputs to be static feedback equivalent to the proposed triangular form, and thus a sufficient condition for an affine input system to be flat. 
	\end{abstract}
	
	
	\end{frontmatter}
	%===============================================================================
	
	%% There are a number of predefined theorem-like environments in
	%% ifacconf.cls:
	%%
	%% \begin{thm} ... \end{thm}            % Theorem
	%% \begin{lem} ... \end{lem}            % Lemma
	%% \begin{claim} ... \end{claim}        % Claim
	%% \begin{conj} ... \end{conj}          % Conjecture
	%% \begin{cor} ... \end{cor}            % Corollary
	%% \begin{fact} ... \end{fact}          % Fact
	%% \begin{hypo} ... \end{hypo}          % Hypothesis
	%% \begin{prop} ... \end{prop}          % Proposition
	%% \begin{crit} ... \end{crit}          % Criterion
	
	\section{Introduction}
		\vspace{-2ex}
		The concept of flatness was introduced in control theory by Fliess, L\'evine, Martin and Rouchon, see \eg \cite{FliessLevineMartinRouchon:1992,FliessLevineMartinRouchon:1995}, and has attracted a lot of interest in the control systems theory community. The flatness property allows an elegant systematic solution of feed-forward and feedback problems, see \eg \cite{FliessLevineMartinRouchon:1995}. Roughly speaking, a nonlinear control system
		\begin{align*}
			\begin{aligned}
				\dot{x}&=f(x,u)
			\end{aligned}
		\end{align*}
		with $\Dim{x}=n$ states and $\Dim{u}=m$ inputs is flat, if there exist $m$ differentially independent functions $y^j=\varphi^j(x,u,u_1,\ldots,u_q)$, $u_k$ denoting the $k$-th time derivative of $u$, such that $x$ and $u$ can be parameterized by $y$ and its time derivatives. In contrast to the static feedback linearization problem, which is completely solved, see \cite{JakubczykRespondek:1980}, \cite{NijmeijervanderSchaft:1990}, there are many open problems concerning flatness. Recent research in the field of flatness can be found in \eg \cite{SchoberlRiegerSchlacher:2010}, \cite{SchlacherSchoberl:2013}, \cite{LiXuSuChu:2013}, \cite{SchoberlSchlacher:2014}, \cite{NicolauRespondek:2014}, \cite{KolarSchoberlSchlacher:2015}, \cite{NicolauRespondek:2017}.
		
		Structurally flat triangular forms are of special interest in the problem of finding flat outputs for nonlinear control systems. Systems that are static feedback linearizable can be transformed into Brunovsky normal form, see \cite{JakubczykRespondek:1980}. Normal forms for systems that become static feedback linearizable after a one-fold prolongation of a suitably chosen control can be found in \cite{NicolauRespondek:2019} (a complete solution of the flatness problem of this class of systems can be found in \cite{NicolauRespondek:2017}). In \cite{BououdenBoutatZhengBarbotKratz:2011}, a structurally flat triangular form for a class of $0$-flat systems is proposed. A structurally flat implicit triangular form for $1$-flat systems is presented in \cite{SchoberlSchlacher:2014}. The flatness problem for two-input driftless systems has been solved in \cite{MartinRouchon:1994}. There, it is shown that a two-input driftless system is flat, if and only if it is static feedback equivalent to the so called chained form (also referred to as Goursat normal form). In \cite{LiXuSuChu:2013} an extension of the chained form to a triangular form, referred to as extended chained form, is considered. In \cite{SilveiraPereiraRouchon:2015} geometric necessary and sufficient conditions for a two-input affine input system (AI-system) to be static feedback equivalent to this extended chained form are derived. Conditions for the case with $m\geq 2$ inputs are provided in \cite{LiNicolauRespondek:2016}, \cite{NicolauLiRespondek:2014} and \cite{Nicolau:2014}.
		
		In this paper, we confine ourselves to AI-systems with two inputs. We present a triangular form based on the extended chained form and derive necessary and sufficient conditions for an AI-system to be static feedback equivalent to this triangular form. This provides a sufficient condition for an AI-system to be flat. The proposed triangular form contains the (extended) chained form as a special case.\vspace{-3ex}
	\section{Notation}
		\vspace{-2ex}
		Let $\mathcal{X}$ be an $n$-dimensional smooth manifold, equipped with local coordinates $x^i$, $i=1,\ldots,n$. Its tangent bundle and cotangent bundle are denoted by $(\mathcal{T}(\mathcal{X}),\tau_\mathcal{X},\mathcal{X})$ and $(\mathcal{T}^\ast(\mathcal{X}),\tau^\ast_\mathcal{X},\mathcal{X})$. For these bundles we have the induced local coordinates $(x^i,\dot{x}^i)$ and $(x^i,\dot{x}_i)$ with respect to the bases $\{\partial_{x^i}\}$ and $\{\D x^i\}$, respectively. We make use of the Einstein summation convention. By $\D\omega$ we denote the exterior derivative of a $p$-form $\omega$. The $k$-fold Lie derivative of a function $\varphi$ along a vector field $v$ is denoted by $\Lie_v^k\varphi$. Let $v,w$ be two vector fields. We denote their Lie bracket by $[v,w]$. Let furthermore $D_1,D_2$ be two distributions. By $[v,D_1]$ we denote the distribution spanned by the Lie bracket of $v$ with all basis vector fields of $D_1$, and by $[D_1,D_2]$ the distribution spanned by the Lie brackets of all possible pairs of basis vector fields of $D_1$ and $D_2$. The $i$-th derived flag of a distribution $D$ is denoted by $D^{(i)}$ and defined by $D^{(0)}=D$ and $D^{(i+1)}=D^{(i)}+[D^{(i)},D^{(i)}]$ for $i\geq 0$. The $i$-th Lie flag of a distribution $D$ is denoted by $D_{(i)}$ and defined by $D_{(0)}=D$ and $D_{(i+1)}=D_{(i)}+[D,D_{(i)}]$ for $i\geq 0$. The involutive closure of $D$ is denoted by $\overline{D}$, it is the smallest involutive distribution containing $D$ and can be determined via the derived flag. We write $\C{D}$ for the Cauchy characteristic distribution of $D$. It is spanned by all vector fields $v$ which belong to $D$ and satisfy $[v,D]\subset D$ and is always involutive. Cauchy characteristics allow us to find a basis for a distribution which is independent of certain coordinates. Since $\C{D}$ is involutive, it can be straightened out such that $\C{D}=\Span{\partial_{x^1},\ldots,\partial_{x^{n_c}}}$, with $n_c=\Dim{\C{D}}$. From $[\C{D},D]\subset D$, it follows that in these coordinates, there exists a basis for $D$ which is independent of the coordinates $(x^1,\ldots,x^{n_c})$. An AI-system
		\vspace{-1.5ex}
		\begin{align}\label{eq:aisystem_intro}
			\begin{aligned}
				\dot{x}&=a(x)+b_j(x)u^j\,,&j&=1,\ldots,m\,,
			\end{aligned}
		\end{align}
		is characterized by the drift vector field $a=a^i(x)\partial_{x^i}$ and the input vector fields $b_j=b_j^i(x)\partial_{x^i}$, $j=1,\ldots,m$, $i=1,\ldots,n$ on the state manifold $\mathcal{X}$. Throughout, we assume that these vector fields and all functions we deal with are smooth. Given two AI-systems, we call them static feedback equivalent, if they are equivalent via a diffeomorphism $\bar{x}=\Phi(x)$ on the state space and an invertible feedback transformation $\bar{u}^j=g^j(x)+m^j_k(x)u^k$. The equivalent system reads
		\begin{align*}
			\begin{aligned}
				\dot{\bar{x}}&=\bar{a}(\bar{x})+\bar{b}_j(\bar{x})\bar{u}^j\,,&j&=1,\ldots,m\,.
			\end{aligned}
		\end{align*}
		\vspace{-4ex}
	\section{Known results}\label{se:knownResults}
	\vspace{-2ex}
		In this section, we summarize some known results from the literature which are of particular importance for characterizing our triangular form. Throughout, we assume all distributions to have locally constant dimension, we consider generic points only. In particular, we call a system static feedback equivalent to a certain normal form, even though the transformation into this form may exhibit singularities. A system is static feedback linearizable, if it is static feedback equivalent to a linear controllable system, in particular to the Brunovsky normal form. The static feedback linearization problem has been solved in \cite{JakubczykRespondek:1980} and \cite{HuntSu:1981}. The geometric necessary and sufficient conditions read as follows. For \eqref{eq:aisystem_intro}, define the distributions $D_{i+1}=D_i+[a,D_i]$, $i\geq 1$, where $D_1=\Span{b_1,\ldots,b_m}$. An $m$-input AI-system \eqref{eq:aisystem_intro} is static feedback linearizable if and only if all the distributions $D_i$, $i\geq 1$ are involutive and $D_{n-1}=\mathcal{T}(\mathcal{X})$.
		 
		The flatness problem for two-input driftless systems %of the form
		\begin{align}\label{eq:driftless}
			\begin{aligned}
				\dot{x}&=b_1(x)u^1+b_2(x)u^2
			\end{aligned}
		\end{align}
		has been solved in \cite{MartinRouchon:1994}. There, it is shown that \eqref{eq:driftless} is flat, if and only if it is static feedback equivalent to the chained form
		\begin{align}\label{eq:chainedform}
			\begin{aligned}
			\dot{x}^1=u^2,~\dot{x}^2=x^3u^2,~\cdots,~\dot{x}^{n-1}=x^nu^2,~\dot{x}^n=u^1\,.
			\end{aligned}
		\end{align}
		The input vector fields of a system in chained form read
		\begin{align}\label{eq:chainedInputvectorfields}
			\begin{aligned}
				b_1&=\partial_{x^n}\,,&b_2&=\partial_{x^1}+x^3\partial_{x^2}+\ldots+x^n\partial_{x^{n-1}}\,.
			\end{aligned}
		\end{align}
		The geometric necessary and sufficient condition for \eqref{eq:driftless} to be static feedback equivalent to the chained form is $\Dim{D^{(i)}}=2+i$, $i=0,\ldots,n-2$ with $D=\Span{b_1,b_2}$. In \cite{Murray:1994} it is shown that locally around a point of the state space at which the additional regularity condition $\Dim{D_{(i)}}=2+i$ on the Lie flag of $D$ holds, the transformation into chained form does not exhibit singularities. A system in chained form is obviously flat with the top variables $(x^1,x^2)$ forming a possible flat output. For a comprehensive analysis of the flatness of systems static feedback equivalent to the chained form, a characterization of all their $x$-flat outputs and their singularities, we refer to \cite{LiRespondek:2012}. 

		In \cite{LiXuSuChu:2013} an extension of the chained form to the triangular form (referred to as extended chained form)
		\begin{align}\label{eq:extendedChainedform}
			\begin{aligned}
				\dot{x}^1&=u^2\\
				\dot{x}^2&=x^3u^2+a^2(x^1,x^2,x^3)\\
				\dot{x}^3&=x^4u^2+a^3(x^1,x^2,x^3,x^4)\\
				&\vdotswithin{=}\\
				\dot{x}^{n-1}&=x^nu^2+a^{n-1}(x^1,\ldots,x^n)\\
				\dot{x}^n&=u^1\,.
			\end{aligned}
		\end{align}
		is considered. In \cite{Nicolau:2014} and \cite{SilveiraPereiraRouchon:2015} geometric necessary and sufficient conditions for an AI-system 
		\vspace{-1.5ex}
		\begin{align*}
			\begin{aligned}
				\dot{x}&=a(x)+b_1(x)u^1+b_2(x)u^2
			\end{aligned}
		\end{align*}
		to be static feedback equivalent to \eqref{eq:extendedChainedform} are provided. Those read as $\Dim{D^{(i)}}=\Dim{D_{(i)}}=2+i$, $i=0,\ldots,n-2$ with $D=\Span{b_1,b_2}$, \ie that for $a(x)=0$, the system is static feedback equivalent to the chained form \eqref{eq:chainedform} and
		\begin{align}\label{eq:compatibility}
			\begin{aligned}
				[a,\mathcal{C}(D^{(i)})]&\subset D^{(i)}\,,&i&=1,\ldots,n-3\,,
			\end{aligned}
		\end{align}
		which assures that the drift is compatible with the chained form. Thus, every flat output of the driftless system obtained by setting $a(x)=0$, is also a flat output of \eqref{eq:extendedChainedform}. For a comprehensive analysis of the flatness of systems static feedback equivalent to the extended chained form, a characterization of their flat outputs and their singularities, we refer to \cite{LiNicolauRespondek:2016} or \cite{Nicolau:2014}.
		
		Transformations to the (extended) chained form can be found in the literature. However, for our purposes, a different procedure for successively transforming a system into the (extended) chained form is beneficial (will be used in the proof of the main theorem) which is described in Appendix \ref{ap:supplements}.
		\vspace{-2ex}
	\section{Flat triangular form based on the extended chained form}
		\vspace{-2ex}
		In the following, we consider the structurally flat triangular form
		\vspace{-1.5ex}
		\begin{align}\label{eq:triangularform}
			\begin{aligned}
				\dot{x}_1&=f_1(x_1,x_2^1,x_2^2)\\
				\dot{x}_2&=f_2(x_1,x_2,x_{3,1}^1,x_{3,2}^1)\\
				\dot{x}_3&=f_3(x_3,u^1,u^2)
			\end{aligned}
		\end{align}
		with the $x_1$-subsystem being in Brunovsky normal form
		\begin{align}\label{eq:subsys1}
			\begin{aligned}
				f_1:\quad\begin{aligned}
					\dot{x}_{1,1}^1&=x_{1,1}^2&\dot{x}_{1,2}^1&=x_{1,2}^2\\
					\dot{x}_{1,1}^2&=x_{1,1}^3&\dot{x}_{1,2}^2&=x_{1,2}^3\\
					&\vdotswithin{=}&&\vdotswithin{=}\\
					\dot{x}_{1,1}^{n_{1,1}}&=x_2^1&\dot{x}_{1,2}^{n_{1,2}}&=x_2^2\,,
				\end{aligned}
			\end{aligned}
		\end{align}
		the $x_2$-subsystem being essentially in extended chained form\vspace{-1.5ex}
		\begin{align}\label{eq:subsys2}
			\begin{aligned}
				f_2:\quad\begin{aligned}
					\dot{x}_2^1&=x_{3,2}^1\\
					\dot{x}_2^2&=x_2^3x_{3,2}^1+a_2^2(x_1,x_2^1,\ldots,x_2^3)\\
					&\vdotswithin{=}\\
					\dot{x}_2^{n_2-1}&=x_2^{n_2}x_{3,2}^1+a_2^{n_2-1}(x_1,x_2)\\
					\dot{x}_2^{n_2}&=x_{3,1}^1
				\end{aligned}
			\end{aligned}
		\end{align}
		and the $x_3$-subsystem again being in Brunovsky normal form\vspace{-1.5ex}
		\begin{align}\label{eq:subsys3}
			\begin{aligned}
				f_3:\quad\begin{aligned}
						\dot{x}_{3,1}^1&=x_{3,1}^2&\dot{x}_{3,2}^1&=x_{3,2}^2\\
					\dot{x}_{3,1}^2&=x_{3,1}^3&\dot{x}_{3,2}^2&=x_{3,2}^3\\
					&\vdotswithin{=}&&\vdotswithin{=}\\
					\dot{x}_{3,1}^{n_3}&=u^1&\dot{x}_{3,2}^{n_3}&=u^2\,.
				\end{aligned}
			\end{aligned}
		\end{align}
		The triangular form \eqref{eq:triangularform} consists of three subsystems. The $x_1$-subsystem is in Brunovsky normal form, it consists of two integrator chains of arbitrary lengths $n_{1,1}\geq 0$ and $n_{1,2}\geq 0$. In total, it consists of $n_1=n_{1,1}+n_{1,2}\geq 0$ states. The $x_2$-subsystem is essentially in extended chained form (the minor difference is that the drift of the $x_2$-subsystem may also depend on the states $x_1$), we assume $n_2\geq 3$. The top variables $x_2^1$ and $x_2^2$ of the $x_2$-subsystem act as inputs for the $x_1$-subsystem. The $x_3$-subsystem is again in Brunovsky normal form, it consists of two integrator chains of equal length $n_3\geq 0$. The top variables $x_{3,1}^1$ and $x_{3,2}^1$ act as inputs for the $x_2$-subsystem (for $n_3=0$, the $x_3$-subsystem does not exist and in the $x_2$-subsystem $x_{3,1}^1$ and $x_{3,2}^1$ are replaced by $u^1$ and $u^2$). In conclusion, the $x_3$-subsystem and the $x_2$-subsystem form an endogenous dynamic feedback for the $x_1$-subsystem. The $x_3$-subsystem in turn is an endogenous dynamic feedback for the $x_2$-subsystem. The total number of states of \eqref{eq:triangularform} is given by $n=n_1+n_2+2n_3$.
		As a motivating example, consider the model of an induction motor
		\begin{align}\label{eq:motor}
			\begin{aligned}
				\dot{\theta}&=\omega\\
				\dot{\omega}&=\mu\psi_dI_q-\tfrac{\tau_L}{J}\\
				\dot{\psi}_d&=-\eta\psi_d+\eta MI_d
			\end{aligned}&\qquad
			\begin{aligned}
				\dot{\rho}&=n_p\omega+\tfrac{\eta MI_q}{\psi_d}\\
				\dot{I}_d&=v_d\\
				\dot{I}_q&=v_q\,,
			\end{aligned}
		\end{align}
		see also \cite{Chiasson:1998}, \cite{NicolauRespondek:2013}, \cite{Schoberl:2014}. This system is not static feedback linearizable, but it is known to be flat with $(\theta,\rho)$ forming a possible flat output. In contrast to a reduced order model discussed in \cite{SilveiraPereiraRouchon:2015}, the system \eqref{eq:motor} is not static feedback equivalent to the extended chained form. However, it is static feedback equivalent to the triangular form \eqref{eq:triangularform} with $n_1=1$, $n_2=3$ and $n_3=1$. Indeed, applying a suitable state and input transformation, \eqref{eq:motor} takes the form
		\vspace{-1ex}
		\begin{align}\label{eq:motorTriangular}
			\begin{aligned}
				&f_1:\quad\begin{aligned}
					\dot{x}_{1,1}^1&=x_2^1
				\end{aligned}\\[1ex]
				&f_2:\quad\begin{aligned}
					\dot{x}_2^1&=x_{3,2}^1\\
					\dot{x}_2^2&=x_2^3x_{3,2}^1+n_p x_2^1+\tfrac{\tau_L}{J}x_2^3\\
					\dot{x}_2^3&=x_{3,1}^1
				\end{aligned}\\[1ex]
				&f_3:\quad\begin{aligned}
					\dot{x}_{3,1}^1&=\bar{u}^1&\dot{x}_{3,2}^1&=\bar{u}^2\,.
				\end{aligned}
			\end{aligned}
		\end{align}
		\begin{rem}
			The system \eqref{eq:motor} becomes static feedback linearizable after a one-fold prolongation of a suitably chosen control. Thus, it is also static feedback equivalent to normal forms presented in \cite{NicolauRespondek:2019}.
		\end{rem}
		\vspace{-2ex}
	\section{Characterization of the triangular form}
		\vspace{-2ex}
		The following theorem provides necessary and sufficient conditions for an AI-system to be static feedback equivalent to the triangular form \eqref{eq:triangularform} and thus provides a sufficient condition for an AI-system to be flat. We again assume all distributions to have locally constant dimension and omit discussing singularities coming along with flat outputs of \eqref{eq:triangularform} or singularities in the problem of transforming a given system into the form \eqref{eq:triangularform}. We consider generic points only, regularity conditions are omitted. Consider a two-input AI-system
		\begin{align}\label{eq:aisystem}
			\begin{aligned}
				\dot{x}&=a(x)+b_1(x)u^1+b_2(x)u^2
			\end{aligned}
		\end{align}
		and let us define the distributions $D_i$, $i=1,\ldots,n_3+1$ where $D_1=\Span{b_1,b_2}$ and $D_{i+1}=D_i+[a,D_i]$, with the smallest integer $n_3$ such that $D_{n_3+1}$ is not involutive. All these distributions are invariant with respect to invertible feedback transformations (this can be shown the same way as in the proof of Proposition 7.1 in \cite{NicolauRespondek:2016}).
		\begin{thm}\label{thm:1}
			The AI-system \eqref{eq:aisystem} is static feedback equivalent to the triangular form \eqref{eq:triangularform} if and only if the following conditions are satisfied:
			\begin{enumerate}[label=(\alph*)]
				\setlength{\itemsep}{5pt}
				\item\label{it:a} $\C{D_{n_3+1}}=D_{n_3}$.
				\item\label{it:b} The derived flags of the non involutive distribution $D_{n_3+1}$ satisfy
				\begin{align*}
					\begin{aligned}
						\Dim{D_{n_3+1}^{(i)}}&=\Dim{D_{n_3+1}}+i\,,\!&i&=1,\ldots,n_2-2,
					\end{aligned}
				\end{align*}
				with the smallest integer $n_2$ such that $D_{n_3+1}^{(n_2-2)}=\overline{D}_{n_3+1}$.
				\item\label{it:c} The drift satisfies the compatibility conditions\footnote{If $\overline{D}_{n_3+1}=\mathcal{T}(\mathcal{X})$, the condition \eqref{eq:coupling} and the items \ref{it:d} and \ref{it:e} have to be omitted. In this case, the system is static feedback equivalent to \eqref{eq:triangularform} with $n_{1,1}=n_{1,2}=0$ if and only if all the other conditions of the theorem are met. (If additionally $n_3=0$ holds, only item \ref{it:b} and condition \eqref{eq:compatibilityThm1} remain, which then match the conditions for static feedback equivalence to the extended chained form.)}
				\begin{align}
					&\!\begin{aligned}\label{eq:compatibilityThm1}
						[a,\C{D_{n_3+1}^{(i)}}]&\subset D_{n_3+1}^{(i)}\,,&i&=1,\ldots,n_2-3\,,
					\end{aligned}\\
					&\begin{aligned}\label{eq:coupling}
						\Dim{\overline{D}_{n_3+1}+[a,D_{n_3+1}^{(n_2-3)}]}&=\Dim{\overline{D}_{n_3+1}}+1\,.\hspace{-2.5ex}
					\end{aligned}
				\end{align}
				\item\label{it:d} The distributions $G_{i+1}$, $i\geq 0$ are involutive, where $G_{i+1}=G_i+[a,G_i]$ and $G_0=\overline{D}_{n_3+1}$.
				\item\label{it:e} $G_s=\mathcal{T}(\mathcal{X})$ holds for some integer $s$.
			\end{enumerate}
		\end{thm}
		\vspace{-1ex}
		All these conditions are easily verifiable and require differentiation and algebraic operations only. A proof of this theorem is provided in Section \ref{subse:proof}. In the following, we outline the meaning of the individual conditions. Consider a system of the form \eqref{eq:triangularform}. The main idea of Theorem \ref{thm:1} is to characterize the individual subsystems separately and have additional conditions which take into account their coupling. The static feedback linearizable subsystem \eqref{eq:subsys3} is characterized by the involutive distributions $D_1,\ldots,D_{n_3}$, we have $D_{n_3}=\Span{\partial_{x_3}}$. Item \ref{it:a} is crucial for the coupling of the $x_2$-subsystem with the $x_3$-subsystem, it guarantees that the $x_2$-subsystem indeed allows an AI representation with respect to its inputs $x_{3,1}^{1}$ and $x_{3,2}^1$. The associated input vector fields of the $x_2$-subsystem are $b_1^c=\partial_{x_2^{n_2}}$ and $b_2^c=\partial_{x_2^1}+x_2^3\partial_{x_2^2}+\ldots+x_2^{n_2}\partial_{x_2^{n_2-1}}$, which are structurally of the form \eqref{eq:chainedInputvectorfields}. We have $D_{n_3+1}=\Span{\partial_{x_3},b_1^c,b_2^c}$ and $\overline{D}_{n_3+1}=\Span{\partial_{x_3},\partial_{x_2}}$. Item \ref{it:b} therefore reflects the fact that the $x_2$-subsystem is essentially in (extended) chained form. Condition \eqref{eq:compatibilityThm1} of item \ref{it:c} assures the compatibility of the drift of the $x_2$-subsystem with its chained structure. Although the drift vector field of \eqref{eq:triangularform} 
		% $a=f_1+x_{3,1}^1b_1^c+x_{3,2}^1b_2^c+a_2+a_3$, where $f_1=x_{1,1}^2\partial_{x_{1,1}^1}+\ldots+x_2^1\partial_{x_{1,1}^{n_{1,1}}}+x_{1,2}^2\partial_{x_{1,2}^1}+\ldots+x_2^2\partial_{x_{1,2}^{n_{1,2}}}$ and $a_3=x_{3,1}^2\partial_{x_{3,1}^1}+\ldots+x_{3,1}^{n_3}\partial_{x_{3,1}^{n_3-1}}+x_{3,2}^2\partial_{x_{3,2}^1}+\ldots+x_{3,2}^{n_3}\partial_{x_{3,2}^{n_3-1}}$, of a system of the form \eqref{eq:triangularform} 
		does not only consist of the drift vector field $a_2=a_2^2\partial_{x_2^2}+\ldots+a_2^{n_2-1}\partial_{x_2^{n_2-1}}$ of the $x_2$-subsystem, there still applies a compatibility condition analogous to \eqref{eq:compatibility} for the extended chained form. Condition \eqref{eq:coupling} of item \ref{it:c} is crucial for the coupling of the $x_2$-subsystem with the $x_3$-subsystem. The items \ref{it:d} and \ref{it:e} characterize the static feedback linearizable $x_3$-subsystem.
		\vspace{-1.5ex}
		\subsection{Determining flat outputs}\label{subse:flatOutputs}
		\vspace{-2ex}
			For determining flat outputs of a system which is static feedback equivalent to \eqref{eq:triangularform}, there is no need to actually transform the system into the form \eqref{eq:triangularform}. Flat outputs are determined directly from the distributions involved in the conditions of Theorem \ref{thm:1}. We have to distinguish between three cases regarding the lengths of the integrator chains of the $x_1$-subsystem in a corresponding representation in the form \eqref{eq:triangularform}. Namely, the cases that 1) both integrator chains have at least a length of one, \ie $n_{1,1},n_{1,2}\geq 1$,\linebreak 2) both integrator chains have length zero, \ie $n_1=0$, 3) one of the integrator chains has length zero, the other length $n_1\geq 1$. We can easily test which case applies. If we have $\Dim{G_1}=\Dim{\overline{D}_{n_3+1}}+2$, both chains have at least length one, if $\overline{D}_{n_3+1}=\mathcal{T}(\mathcal{X})$, both have length zero, if $\Dim{G_1}=\Dim{\overline{D}_{n_3+1}}+1$, one chain has length zero. In the following, we discuss these three cases in more detail.
			\vspace{-1.2ex}
			\subsubsection{Case 1:} If $n_{1,1},n_{1,2}\geq 1$, \ie both integrator chains of the $x_1$-subsystem have at least length one, flat outputs of \eqref{eq:triangularform} are determined from the sequence of involutive distributions of Theorem \ref{thm:1} item \ref{it:d} in the same way as linearizing outputs are determined from the sequence of involutive distributions involved in the test for static feedback linearizability (see \eg \cite{JakubczykRespondek:1980}, \cite{NijmeijervanderSchaft:1990}). 
			\vspace{-1.2ex}
			\subsubsection{Case 2:} If both chains have length zero (no $x_1$-subsystem exists), the problem of finding flat outputs of \eqref{eq:triangularform} is in fact the same as finding flat outputs of a system which is static feedback equivalent to the chained form. This problem is addressed in \cite{LiRespondek:2012}. In this case, flat outputs are all pairs of functions $(\varphi^1,\varphi^2)$, which meet $L=(\Span{\D\varphi^1,\D\varphi^2})^\perp\subset D_{n_3+1}^{(n_2-3)}$, with $D_{n_3+1}^{(n_2-3)}$ of Theorem \ref{thm:1} item \ref{it:b}. In \cite{LiRespondek:2012}, Theorem 2.10, a method for constructing such a distribution $L$ is provided. The distribution $L$ is not unique, one has to choose one function $\varphi^1$ whose differential $\D\varphi^1\neq 0$ annihilates $\mathcal{C}(D_{n_3+1}^{(n_2-3)})$. Once such a function has been chosen, the distribution $L$ is uniquely determined by this choice and can be calculated, and in turn a possible second function $\varphi^2$, which together with $\varphi^1$ forms a possible flat output, can be calculated. Equivalent to the method for determining $L$ provided in \cite{LiRespondek:2012}, once a function $\varphi^1$ whose differential annihilates $\mathcal{C}(D_{n_3+1}^{(n_2-3)})$ has been chosen, the annihilator of $L$ can also be calculated via $L^\perp=(D_{n_3+1}^{(n_2-3)})^\perp+\Span{\D\varphi^1}$.
			\vspace{-1.2ex}
			\subsubsection{Case 3:} Here, the $x_1$-subsystem determines one component $\varphi^1$ of a flat output. This function is obtained by integrating $G_{s-1}^\perp$, \ie by finding a function $\varphi^1$ such that $\Span{\D\varphi^1}=G_{s-1}^\perp$, with $G_{s-1}$ of Theorem \ref{thm:1}, item \ref{it:d}. A possible second component $\varphi^2$ is obtained by integrating the integrable codistribution $L^\perp=(D_{n_3+1}^{(n_2-3)})^\perp+\Span{\D\Lie_a^s\varphi^1}$. A function $\varphi^2$ whose differential $\D\varphi^2$ together with $\{\D\varphi^1,\D\Lie_a\varphi^1,\ldots,\D\Lie_a^s\varphi^1\}$ spans the codistribution $L^\perp$, is a possible second component. In this case, the distribution $L$ is uniquely determined by the function $\Lie_a^s\varphi^1$, which is imposed by the $x_1$-subsystem.
			\vspace{-1ex}
		\subsection{Remarks and limitations}\label{subse:remarks}
		\vspace{-2ex}
			As already mentioned, the triangular form \eqref{eq:triangularform} contains the (extended) chained form as a special case. A system in chained form \eqref{eq:chainedform} can be rendered static feedback linearizable by a $(n-2)$-fold prolongation of the control $u^2$. The same applies to a system in extended chained form \eqref{eq:extendedChainedform} and analogously, a system in the triangular form \eqref{eq:triangularform} becomes static feedback linearizable by a $(n_2-2)$-fold prolongation of the control $u^2$. Therefore, the case $n_2=3$ is covered by \cite{NicolauRespondek:2013}. The case $n_2=4$ is covered by \cite{NicolauRespondek:2016-2} if additionally $\overline{D}_{n_3+1}\neq\mathcal{T}(\mathcal{X})$ holds, \ie if $n_1\geq 1$. Our geometric characterization of \eqref{eq:triangularform} provided by Theorem \ref{thm:1} is not subject to restrictions on $n_1$ or $n_2$. An obvious restriction of the triangular form \eqref{eq:triangularform} is however that the integrator chains of the $x_3$-subsystem are supposed to have equal lengths. Even though the case where the lengths differ by at most one seems to be tractable, we do not discuss it in this contribution.	
			\vspace{-1.5ex}	
		\subsection{Proof of Theorem \ref{thm:1}}\label{subse:proof}
			\vspace{-1.5ex}
			The following two lemmas are of particular importance for the sufficiency part of the proof.
			\begin{lem}\label{lem:characteristic}
				Let $D$ be a distribution. Every characteristic vector field of $D$, \ie every vector field $c\in\C{D}$ is also characteristic for its derived flag $D^{(1)}$, \ie $\C{D}\subset\C{D^{(1)}}$.
			\end{lem}
			\vspace{-1ex}
			An immediate consequence of Lemma \ref{lem:characteristic} is that the Cauchy characteristic distributions $\C{D^{(i)}}$, $i\geq 0$ form the sequence of nested involutive distributions $\C{D}\subset\C{D^{(1)}}\subset\C{D^{(2)}}\subset\ldots$
			\begin{lem}\label{lem:cartan}
				If a $d$-dimensional distribution $D$ satisfies $\Dim{\C{D}}=d-2$ and $\Dim{D^{(i)}}=d+i$, $i=1,\ldots,l$ with $l$ such that $D^{(l)}=\overline{D}$, then the Cauchy characteristics $\C{D^{(i)}}$ satisfy $\Dim{\C{D^{(i)}}}=d-2+i$ and $\C{D^{(i)}}\subset D^{(i-1)}$, $i=1,\ldots,l-1$.
			\end{lem}
			\vspace{-1ex}
			Lemma \ref{lem:cartan} is based on a similar one in \cite{Cartan:1914}, see also \cite{MartinRouchon:1994}, Lemma 2.
			
			{\it Necessity}. We have to show that a system of the form \eqref{eq:triangularform} fulfills the conditions of Theorem \ref{thm:1}. The drift vector field of \eqref{eq:triangularform} reads $a=f_1+x_{3,1}^1b_1^c+x_{3,2}^1b_2^c+a_2+a_3$ where $f_1=x_{1,1}^2\partial_{x_{1,1}^1}+\ldots+x_2^1\partial_{x_{1,1}^{n_{1,1}}}+x_{1,2}^2\partial_{x_{1,2}^1}+\ldots+x_2^2\partial_{x_{1,2}^{n_{1,2}}}$, $b_1^c=\partial_{x_2^{n_2}}$, $b_2^c=\partial_{x_2^1}+x_2^3\partial_{x_2^2}+\ldots+x_2^{n_2}\partial_{x_2^{n_2-1}}$, $a_2=a_2^2\partial_{x_2^2}+\ldots+a_2^{n_2-1}\partial_{x_2^{n_2-1}}$ and $a_3=x_{3,1}^2\partial_{x_{3,1}^1}+\ldots+x_{3,1}^{n_3}\partial_{x_{3,1}^{n_3-1}}+x_{3,2}^2\partial_{x_{3,2}^1}+\ldots+x_{3,2}^{n_3}\partial_{x_{3,2}^{n_3-1}}$. The input vector fields of \eqref{eq:triangularform} are given by $b_1=\partial_{x_{3,1}^{n_3}}$ and $b_2=\partial_{x_{3,2}^{n_3}}$. The distributions defined right before Theorem \ref{thm:1} are thus given by
			\begin{align*}
				\begin{aligned}
					D_1&=\Span{\partial_{x_{3,1}^{n_3}},\partial_{x_{3,2}^{n_3}}}\\
					D_2&=\Span{\partial_{x_{3,1}^{n_3}},\partial_{x_{3,2}^{n_3}},\partial_{x_{3,1}^{n_3-1}},\partial_{x_{3,2}^{n_3-1}}}\\
					&\vdotswithin{=}\\
					D_{n_3}&=\Span{\underbrace{\partial_{x_{3,1}^{n_3}},\partial_{x_{3,2}^{n_3}},\ldots,\partial_{x_{3,1}^{1}},\partial_{x_{3,2}^1}}_{\partial_{x_3}}}\\
					D_{n_3+1}&=\Span{\partial_{x_3},b_1^c,b_2^c}\,.
				\end{aligned}
			\end{align*}
			\subsubsection{Item \ref{it:a}:} Since the input vector fields $b_1^c$, $b_2^c$ of the $x_2$-subsystem are independent of the states $x_3$, $\C{D_{n_3+1}}=D_{n_3}$ indeed holds.
			\subsubsection{Item \ref{it:b}:} The derived flags of the non involutive distribution $D_{n_3+1}$ are given by
			\begin{align*}
				\begin{aligned}
					D_{n_3+1}^{(i)}&=\Span{\partial_{x_3},\partial_{x_2^1}+x_2^3\partial_{x_2^2}+\ldots+x_2^{n_2-i}\partial_{x_2^{n_2-i-1}},\\
						&\hspace{6em}\partial_{x_2^{n_2}},\ldots,\partial_{x_2^{n_2-i}}}\,,~i=0,\ldots,n_2-2\,,
				\end{aligned}
			\end{align*}
			where $D_{n_3+1}^{(n_2-2)}=\overline{D}_{n_3+1}=\Span{\partial_{x_3},\partial_{x_2}}$ and thus, the condition $\Dim{D_{n_3+1}^{(i)}}=\Dim{D_{n_3+1}}+i$, $i=1,\ldots,n_2-2$ is met. Since $b_1^c$ and $b_2^c$ are independent of $x_3$, the derived flags are solely determined by $b_1^c$ and $b_2^c$, which are structurally of the form \eqref{eq:chainedInputvectorfields}. 
			\vspace{-1.5ex}
			\subsubsection{Item \ref{it:c}:} This item consists of two conditions taking into account the drift of the system. The condition \eqref{eq:compatibilityThm1} corresponds to the compatibility condition \eqref{eq:compatibility} of the extended chained form \eqref{eq:extendedChainedform}. The condition \eqref{eq:coupling} is crucial for the coupling of $x_2$-subsystem with the $x_3$-subsystem. The Cauchy characteristic distributions of the derived flags $D_{n_3+1}^{(i)}$ are given by
			\begin{align*}
				\begin{aligned}
					\C{D_{n_3+1}^{(i)}}=\Span{\partial_{x_3},\partial_{x_2^{n_2}},\ldots,\partial_{x_2^{n_2-i+1}}},\,
				\end{aligned}
			\end{align*}
			for $i=1,\ldots,n_2-3$. That the compatibility condition \eqref{eq:compatibilityThm1}, \ie $[a,\C{D_2^{(i)}}]\subset D_2^{(i)}$, $i=1,\ldots,n_2-3$ is indeed met, follows from
%			\footnote{Evaluated for \eg $i=1$, we obtain
%				\begin{align*}
%					\begin{aligned}
%						[\partial_{x_2^{n_2}},a]&=(x_{3,2}^1+\partial_{x_2^{n_2}}a^{n_2-1})\partial_{x_2^{n_2-1}}\Mod\C{\Delta_1^{(1)}}\,.
%					\end{aligned}
%				\end{align*}
%				Because of $\Delta_1^{(1)}=\Span{\partial_{x_3},\partial_{x_2^1}+x_2^3\partial_{x_2^2}+\ldots+x_2^{n_2-1}\partial_{x_2^{n_2-2}},\partial_{x_2^{n_2}},\partial_{x_2^{n_2-1}}}$ and $\C{\Delta_1^{(1)}}=\Span{\partial_{x_3},\partial_{x_2^{n_2}}}$, it indeed follows that $[a,\C{\Delta_1^{(1)}}]\subset\Delta_1^{(1)}$.}
			\begin{align*}
				\begin{aligned}
					[\partial_{x_2^{n_2-i+1}},a]&=(x_{3,2}^1+\partial_{x_2^{n_2-i+1}}a^{n_2-i})\partial_{x_2^{n_2-i}}\Mod\C{\Delta_1^{(i)}}
				\end{aligned}
			\end{align*}
			for $i=1,\ldots,n_2-3$. The condition \eqref{eq:coupling} follows from
			\begin{align*}
				\begin{aligned}
					\overline{D}_{n_3+1}+&[a,D_{n_3+1}^{(n_2-3)}]=\\
					&~~~~\overline{D}_{n_3+1}+\Span{\underbrace{[a,\partial_{x_2^{n_2}}],\ldots,[a,\partial_{x_2^{3}}]}_{\in\overline{D}_2},\underbrace{[a,b^c_2]}_{\notin\overline{D}_2}}\,,
				\end{aligned}
			\end{align*}
			again with $b^c_2=\partial_{x_2^{1}}+x_2^3\partial_{x_2^{2}}+x_2^4\partial_{x_2^{3}}+\ldots+x_2^{n_2}\partial_{x_2^{n_2-1}}$.
			\vspace{-1.5ex}
			\subsubsection{Item \ref{it:d} and \ref{it:e}:} The distributions $G_i$ follow as
			\begin{align*}
				\begin{aligned}
					G_0&=\overline{D}_{n_3+1}=\Span{\partial_{x_3},\partial_{x_2}}\\
					G_i&=\Span{\partial_{x_3},\partial_{x_2},\partial_{x_{1,1}^{n_{1,1}}},\ldots,\partial_{x_{1,1}^{n_{1,1}-i+1}},\\
						&~~~~~~~~~~~~~~~~~~\partial_{x_{1,2}^{n_{1,2}}},\ldots,\partial_{x_{1,2}^{n_{1,2}-i+1}}}\,,~i\geq 1\,,
				\end{aligned}
			\end{align*}
			where any $\partial_{x_{1,j}^k}$, $j=1,2$ with $k\leq 0$ has to be omitted. These distributions are involutive and we have $G_s=\mathcal{T}(\mathcal{X})$ for $s=\max\{n_{1,1},n_{1,2}\}$.
						
			{\it Sufficiency}. We have to show that an AI-system which meets the conditions of Theorem \ref{thm:1} can be transformed into the triangular form \eqref{eq:triangularform}. Due to space limitations, we only give a sketch of the proof, which consists of several steps. \mbox{The steps are demonstrated on Example 1 in Section \ref{se:examples}.}
			
			Consider a system which meets the conditions of Theorem \ref{thm:1}. Due to Lemma \ref{lem:characteristic}, the Cauchy characteristics of the derived flags of $D_{n_3+1}$ form the sequence of nested involutive distributions $\C{D_{n_3+1}^{(1)}}\subset\ldots\subset\C{D_{n_3+1}^{(n_2-3)}}$. We thus have the following sequence of nested involutive distributions
			\begin{align}\label{eq:involutiveSequenceProof}
				\begin{aligned}
					&\!\! D_1\subset\ldots\subset D_{n_3}\subset\C{D_{n_3+1}^{(1)}}\subset\ldots\subset\\
					&~~~~~\C{D_{n_3+1}^{(n_2-3)}}\subset
					\overline{D}_{n_3+1}\subset G_1\subset\ldots\subset G_s=\mathcal{T}(\mathcal{X})\,.\!\!
				\end{aligned}
			\end{align}
			\begin{rem}
				Regarding the dimensions of these distributions, note that $\C{D_{n_3+1}}=D_{n_3}$ implies that $\Dim{D_i}=2i$, $i=1,\ldots,n_3+1$. Because of item \ref{it:a} and \ref{it:b}, the distribution $D_{n_3+1}$ meets the conditions of Lemma \ref{lem:cartan} and thus, $\Dim{\C{D_{n_3+1}^{(i)}}}=2n_3+i$. By the definition of $n_2$, we have $\Dim{\overline{D}_{n_3+1}}=2n_3+n_2$. For the distributions $G_i$, we have either $\Dim{G_1}=\Dim{\overline{D}_{n_3+1}}+2$ and $1\leq\Dim{G_{i+1}}-\Dim{G_i}\leq2$ or $\Dim{G_1}=\Dim{\overline{D}_{n_3+1}}+1$ and $\Dim{G_{i+1}}=\Dim{G_i}+1$ for $i=1,\ldots,s-1$.
			\end{rem}
			The transformation of \eqref{eq:aisystem} into the form \eqref{eq:triangularform} is done in the following six steps.
			\vspace{-1.5ex}
			\subsubsection{Step 1:} Straighten out all the distributions \eqref{eq:involutiveSequenceProof} simultaneously. In such coordinates, the system \eqref{eq:aisystem} takes the form
			\vspace{-1.5ex}
			\begin{align}\label{eq:step1}
				\begin{aligned}
					\dot{x}_1&=f_1(x_1,x_2^1,x_2^2,x_2^3)\\					\dot{x}_2&=f_2(x_1,x_2,x_3^1,x_3^2)\\
					\dot{x}_3&=f_3(x_1,x_2,x_3,u^1,u^2)%a_3(x)+b_{3,1}(x)u^1+b_{3,2}(x)u^2
					\,,
				\end{aligned}
			\end{align}
			with $\Rank{\partial_{(x_2^1,x_2^2,x_2^3)}f_1}\leq 2$ and $\Rank{\partial_{(x_3^1,x_3^2)}f_2}=2$, \ie the system is already decomposed into three subsystems. The $x_1$-subsystem and the $x_3$-subsystem are already in a triangular form, known from the static feedback linearization problem (see \eg \cite{NijmeijervanderSchaft:1990}), the inputs $u^1$, $u^2$ of course occur affine in $f_3$. That $f_1$ is independent of the states $(x_1,x_2^4,\ldots,x_2^{n_2})$ is implied by condition \eqref{eq:compatibilityThm1} of item \ref{it:c}, evaluated for $i=n_2-3$, \ie $[a,\C{D_{n_3+1}^{(n_2-3)}}]\subset D_{n_3+1}^{(n_2-3)}$. 
			\vspace{-1ex}
			\subsubsection{Step 2:} Transform the $x_1$-subsystem into Brunovsky normal form by successively introducing new coordinates from top to bottom. In the prior to last step, the system reads
			\vspace{-1ex}
			\begin{align*}
				\begin{aligned}
					f_1:\quad&\begin{aligned}
						\dot{x}_{1,1}^1&=x_{1,1}^2&\dot{x}_{1,2}^1&=x_{1,2}^2\\
						\dot{x}_{1,1}^2&=x_{1,1}^3&\dot{x}_{1,2}^2&=x_{1,2}^3\\
						&\vdotswithin{=}&&\vdotswithin{=}\\
						\dot{x}_{1,1}^{n_{1,1}}&=\varphi^1(\bar{x}_1,x_2^1,x_2^2,x_2^3)&\dot{x}_{1,2}^{n_{1,2}}&=\varphi^2(\bar{x}_1,x_2^1,x_2^2,x_2^3)
					\end{aligned}\\[1ex]
					&\quad\begin{aligned}
						\dot{x}_2&=\bar{f}_2(\bar{x}_1,x_2,x_3^1,x_3^2,x_3^3)\\
						\dot{x}_3&=\bar{f}_3(\bar{x}_1,x_2,x_3,u^1,u^2)\,,
					\end{aligned}
				\end{aligned}
			\end{align*}
			with $\bar{x}_1=(x_{1,1}^1,\ldots,x_{1,1}^{n_{1,1}},x_{1,2}^1,\ldots,x_{1,2}^{n_{1,2}})$. The functions $\varphi^j(\bar{x}_1,x_2^1,x_2^2,x_2^3)$, $j=1,2$ of the $x_1$-subsystem determine the desired top variables for the $x_2$-subsystem. These functions meet $\D\bar{x}_1\wedge\D\varphi^1\wedge\D\varphi^2\neq 0$. For the system to be static feedback equivalent to the triangular form \eqref{eq:triangularform}, the top variables of the $x_2$-subsystem have to form a flat output of the $x_2$-subsystem which is compatible with its (extended) chained form. It can be shown that item \ref{it:c} implies that the functions $\varphi^j(\bar{x}_1,x_2^1,x_2^2,x_2^3)$ satisfy $L=(\Span{\D \bar{x}_1,\D\varphi^1,\D\varphi^2})^\perp\subset D_{n_3+1}^{(n_2-3)}$, which in turn implies that they indeed form a flat output compatible with the (extended) chained form\footnote{This case corresponds to Case 1 in Section \ref{subse:flatOutputs}. In Case 2, there does not exist an $x_1$-subsystem which determines the desired top variables for the $x_2$-subsystem, any two functions $\varphi^1$, $\varphi^2$ which form a compatible flat output of the $x_2$-subsystem can thus be chosen as the top variables. In Case 3, \ie in case that the $x_1$-subsystem consists only of one integrator chain, the $x_1$-subsystem determines only one function $\varphi^1$. In this case, there always exists a second independent function $\varphi^2$ which together with $\varphi^1$ forms a compatible flat output of the $x_2$-subsystem.}.
			\vspace{-1ex}
			\subsubsection{Step 3:} Introduce the functions $\varphi^1,\varphi^2$ as the top variables of the $x_2$-subsystem, \ie apply the state transformation $\bar{x}_2^j=\varphi^j$, $j=1,2$, with all the other coordinates left unchanged. This completes the transformation of the $x_1$-subsystem to Brunovsky normal form.
			
			In the following Steps 4 and 5, the $x_2$-subsystem is transformed into (extended) chained form following essentially the procedure described in Appendix \ref{ap:supplements}.
			\subsubsection{Step 4:} Normalize the first equation of the $x_2$-subsystem such that $\dot{\bar{x}}_2^1=x_{3,2}^1$. Based on item \ref{it:a} and Lemma \ref{lem:cartan} applied to $D_{n_3+1}$, it can be shown that 
			after this normalization, the $x_2$-subsystem reads
			\begin{align*}
				f_2:\hspace{-1.5ex}\begin{aligned}
					\dot{\bar{x}}_2^1&=x_{3,2}^1\\
					\dot{\bar{x}}_2^2&=b_2^2(\bar{x}_1,\bar{x}_2^1,\bar{x}_2^2,x_2^3)x_{3,2}^1+a_2^2(\bar{x}_1,\bar{x}_2^1,\bar{x}_2^2,x_2^3)\\
					\dot{x}_2^3&=b_2^3(\bar{x}_1,\bar{x}_2^1,\bar{x}_2^2,x_2^3,x_2^4)x_{3,2}^1+a_2^3(\bar{x}_1,\bar{x}_2^1,\bar{x}_2^2,x_2^3,x_2^4)\\
					&\vdotswithin{=}\\
					\dot{x}_2^{n_2-1}&=b_2^{n_2-1}(\bar{x}_1,\bar{x}_2)x_{3,2}^1+a_2^{n_2-1}(\bar{x}_1,\bar{x}_2)\\
					\dot{x}_2^{n_2}&=g(\bar{x}_1,\bar{x}_2,x_3^1,x_{3,2}^1)\,.
				\end{aligned}
			\end{align*}
			The triangular dependence of the functions $a_2^i$ on the states of the $x_2$-subsystem follows from evaluating the condition \eqref{eq:compatibilityThm1}, \ie $[a,\C{D_{n_3+1}^{(i)}}]\subset D_{n_3+1}^{(i)}$, $i=1,\ldots,n_2-3$, in these coordinates.
			\vspace{-1ex}
			\subsubsection{Step 5:} Successively introduce the functions $b_2^i$ as new states of the $x_2$-subsystem and complete the transformation of the $x_2$-subsystem into extended chained form by subsequently normalizing its last equation such that $\dot{\bar{x}}_2^{n_2}=x_{3,1}^1$.
			\vspace{-1ex}
			\subsubsection{Step 6:} Transform the $x_3$-subsystem into Brunovsky normal form by successively introducing new coordinates from top to bottom and applying a suitable static feedback.
			\vspace{-3ex}
	\section{Examples}\label{se:examples}
		\vspace{-2ex}
	  	\subsubsection{Example 1.}\label{ex:academic}
	  		Based on the following academic example, we demonstrate the transformation into the triangular form \eqref{eq:triangularform} by following the six steps of the sufficiency part of the proof of Theorem \ref{thm:1}. Consider the system
			\begin{align}\label{eq:academic}
				\begin{aligned}
					\dot{x}^1&=x^4+1&\dot{x}^5&=x^4(x^7-x^8)\\
					\dot{x}^2&=x^3x^4-x^5&\dot{x}^6&=x^7\\
					\dot{x}^3&=x^7-x^8&\dot{x}^7&=u^1\\
					\dot{x}^4&=x^6(x^7-x^8+x^1)&\dot{x}^8&=u^2\,.
				\end{aligned}
			\end{align}
			The input vector fields are given by $b_1=\partial_{x^7}$ and $b_2=\partial_{x^8}$, the drift is given by
			\begin{align*}
				\begin{aligned}
					a&=(x^4+1)\partial_{x^1}+(x^3x^4-x^5)\partial_{x^2}+(x^7-x^8)\partial_{x^3}+\\
					&\qquad~~ x^6(x^7-x^8+x^1)\partial_{x^4}+x^4(x^7-x^8)\partial_{x^5}+x^7\partial_{x^6}\,.
				\end{aligned}
			\end{align*}
			The distribution $D_1=\Span{b_1,b_2}$ is involutive, the distribution
			\begin{align*}
				\begin{aligned}
					D_2&=D_1+[a,D_1]\\
					&=\Span{\partial_{x^8},\partial_{x^7},\partial_{x^6},\partial_{x^3}+x^6\partial_{x^4}+x^4\partial_{x^5}}
				\end{aligned}
			\end{align*}
			is not involutive, so we have $n_3=1$. We have $\C{D_2}=D_1$, thus, item \ref{it:a} is met. The derived flags of $D_2$ are given by
			\begin{align*}
				\begin{aligned}
					D_2^{(1)}&=\Span{\partial_{x^8},\partial_{x^7},\partial_{x^6},\partial_{x^4},\partial_{x^3}+x^4\partial_{x^5}}\\
					D_2^{(2)}&=\Span{\partial_{x^8},\ldots,\partial_{x^3}}=\overline{D}_2\,.
				\end{aligned}
			\end{align*}
			Thus, item \ref{it:b} is met with $n_2=4$. The drift vector field $a$ satisfies the compatibility condition \eqref{eq:compatibilityThm1} of item \ref{it:c}, which we do not present in detail here. Also the condition \eqref{eq:coupling} of item \ref{it:c} is met. Evaluating item \ref{it:d} yields the distribution $G_1=\overline{D}_2+[a,\overline{D}_2]=\mathcal{T}(\mathcal{X})$, \ie item \ref{it:d} is met and item \ref{it:e} is met with $s=1$. Therefore, according to Theorem \ref{thm:1}, the system is static feedback equivalent to the triangular form \eqref{eq:triangularform}. Since $\Dim{G_1}=\Dim{\overline{D}_2}+2$, the $x_1$-subsystem in a corresponding triangular from \eqref{eq:triangularform} consists of two integrator chains, since $s=1$, both of these integrator chains have length one. According to Section \ref{subse:flatOutputs}, Case 1, flat outputs of \eqref{eq:academic} which are compatible with the triangular form \eqref{eq:triangularform} are all pairs of functions $(\varphi^1,\varphi^2)$ which satisfy $\Span{\D\varphi^1,\D\varphi^2}=G_0^\perp=(\overline{D}_2)^\perp$. We have $(\overline{D}_2)^\perp=\Span{\D x^1,\D x^2}$. A possible flat output is thus \eg $\varphi^1=x^1$, $\varphi^2=x^2$. In the following, we transform \eqref{eq:academic} into the triangular form \eqref{eq:triangularform} by following the six steps of the sufficiency part of the proof of Theorem \ref{thm:1}.
			\vspace{-1ex}
			\subsubsection{Step 1:} In this example, the distributions
			\begin{align*}
				\begin{aligned}
					D_1\subset\C{D_2^{(1)}}\subset\overline{D}_2\subset G_1=\mathcal{T}(\mathcal{X})\,,
				\end{aligned}
			\end{align*}
			which correspond to the sequence \eqref{eq:involutiveSequenceProof}, are already straightened out, \eqref{eq:academic} is therefore structurally already in the form \eqref{eq:step1}. Indeed, renaming the states according to $x_1^1=x^1$, $x_1^2=x^2$, $x_2^1=x^3$, $x_2^2=x^4$, $x_2^3=x^5$, $x_2^4=x^6$, $x_3^1=x^7$ and $x_3^2=x^8$, we obtain
			\begin{align*}
				&f_1:~~\begin{aligned}
						\dot{x}_1^1&=x_2^2+1\\
						\dot{x}_1^2&=x_2^1x_2^2-x_2^3
					\end{aligned}&f_2:~~\begin{aligned}
						\dot{x}_2^1&=x_3^1-x_3^2\\
						\dot{x}_2^2&=x_2^4(x_3^1-x_3^2+x_1^1)\\
						\dot{x}_2^3&=x_2^2(x_3^1-x_3^2)\\
						\dot{x}_2^4&=x_3^1
					\end{aligned}\\
				&f_3:~~\begin{aligned}
					\dot{x}_3^1&=u^1\\
					\dot{x}_3^2&=u^2\,.
				\end{aligned}
			\end{align*}
			\subsubsection{Step 2 and 3:} The $x_1$-subsystem is already in Brunovsky normal form, except for a normalization of the "inputs" of the integrators. To be consistent with the notation in the proof, we only have to rename the states of the $x_1$-subsystem according to $x_{1,1}^1=x_1^1$ and $x_{1,2}^1=x_1^2$. To normalize the "inputs" of the integrators, we introduce $\bar{x}_2^1=x_2^2+1$ and $\bar{x}_2^2=x_2^1x_2^2-x_2^3$ and obtain
			\vspace{0.5ex}
			\begin{align*}
				&f_1:~~\begin{aligned}
					\dot{x}_{1,1}^1&=\bar{x}_2^1\\
					\dot{x}_{1,2}^1&=\bar{x}_2^2
				\end{aligned}&&&f_2:~~\begin{aligned}
					\dot{\bar{x}}_2^1&=x_2^4(x_3^1-x_3^2+x_{1,1}^1)\\
					\dot{\bar{x}}_2^2&=\tfrac{x_2^4(\bar{x}_2^2+x_2^3)}{\bar{x}_2^1-1}(x_3^1-x_3^2+x_{1,1}^1)\\
					\dot{x}_2^3&=(\bar{x}_2^1-1)(x_3^1-x_3^2)\\
					\dot{x}_2^4&=x_3^1\,.
				\end{aligned}
			\end{align*}
			\vspace{-1.5ex}
			\subsubsection{Step 4:} Normalizing the first equation of the $x_2$-subsystem, \ie introducing $x_{3,2}^1=x_2^4(x_3^1-x_3^2+x_{1,1}^1)$ results in
			\vspace{-1.5ex}
			\begin{align*}
				&f_2:\quad\begin{aligned}
					\dot{\bar{x}}_2^1&=x_{3,2}^1\\
					\dot{\bar{x}}_2^2&=\tfrac{\bar{x}_2^2+x_2^3}{\bar{x}_2^1-1}x_{3,2}^1\\
					\dot{x}_2^3&=\tfrac{\bar{x}_2^1-1}{x_2^4}x_{3,2}^1+x_{1,1}^1(1-\bar{x}_2^1)\\
					\dot{x}_2^4&=x_3^1
				\end{aligned}\\[2ex]
				&\!\!\!f_3:\quad\begin{aligned}
					\dot{x}_3^1&=u^1\\
					\dot{x}_{3,2}^1&=%g(\bar{x}_2^1,x_2^4,x_3^1,x_{3,2}^1,u^1,u^2)
					\tfrac{x_3^1x_{3,2}^1}{x_2^4}+x_2^4(u^1-u^2+\bar{x}_2^1)\,.
				\end{aligned}
			\end{align*}
			\vspace{-1ex}
			\subsubsection{Step 5:} Next, we successively introduce the components of the "input" vector field belonging to the "input" $x_{3,2}^1$ of the $x_2$-subsystem as new coordinates. We start with $\bar{x}_2^3=\tfrac{\bar{x}_2^2+x_2^3}{\bar{x}_2^1-1}$. After one more such steps, the $x_2$-subsystem reads
			\vspace{-1.5ex}
			\begin{align*}
				&f_2:\quad\begin{aligned}
					\dot{\bar{x}}_2^1&=x_{3,2}^1\\
					\dot{\bar{x}}_2^2&=\bar{x}_2^3x_{3,2}^1\\
					\dot{\bar{x}}_2^3&=\bar{x}_2^4x_{3,2}^1-x_{1,1}^1\\
					\dot{\bar{x}}_2^4&=-(\bar{x}_2^4)^2x_3^1
				\end{aligned}
			\end{align*}
			and we complete the transformation of the $x_2$-subsystem into extended chained form by normalizing its last equation, \ie by introducing $x_{3,1}^1=-(\bar{x}_2^4)^2x_3^1$, resulting in
			\begin{align*}
				&f_2:~\begin{aligned}
					\dot{\bar{x}}_2^1&=x_{3,2}^1\\
					\dot{\bar{x}}_2^2&=\bar{x}_2^3x_{3,2}^1\\
					\dot{\bar{x}}_2^3&=\bar{x}_2^4x_{3,2}^1-x_{1,1}^1\\
					\dot{\bar{x}}_2^4&=x_{3,1}^1
				\end{aligned}&f_3:~\begin{aligned}
				\dot{x}_{3,1}^1&=\tfrac{2(x_{3,1}^1)^2}{\bar{x}_2^4}-(\bar{x}_2^4)^2u^1\\
				\dot{x}_{3,2}^1&=\tfrac{\bar{x}_2^1-x_{3,1}^1x_{3,2}^1+u^1-u^2}{\bar{x}_2^4}
				\end{aligned}
			\end{align*}
			\vspace{-2ex}
			\subsubsection{Step 6:} The last step is to transform the $x_3$-subsystem into Brunovsky normal form. This is accomplished by applying a suitable static feedback such that $\dot{x}_{3,1}^1=\bar{u}^1$ and $\dot{x}_{3,2}^1=\bar{u}^2$.
		\subsubsection{Example 2.}\label{ex:motor}
			Consider again our motivating example, the model of an induction motor \eqref{eq:motor}. The input vector fields are given by $b_1=\partial_{I_d}$ and $b_2=\partial_{I_q}$, the drift is given by
			\begin{align*}
				\begin{aligned}
					a&=\omega\partial_{\theta}+(\mu\psi_dI_q-\tfrac{\tau_L}{J})\partial_{\omega}+(-\eta\psi_d+\eta MI_d)\partial_{\psi_d}+\\
					&\quad\quad\quad (n_p\omega+\tfrac{\eta MI_q}{\psi_d})\partial_{\rho}\,.
				\end{aligned}
			\end{align*}
			In the following, we apply Theorem \ref{thm:1} to show that this system is indeed static feedback equivalent to the triangular form \eqref{eq:triangularform}. The distribution $D_1=\Span{b_1,b_2}$ is involutive, the distribution
			\begin{align*}
				\begin{aligned}
					D_2&=D_1+[a,D_1]\\
					&=\Span{\partial_{I_d},\partial_{I_q},\partial_{\psi_d},(\psi_d)^2\mu\partial_{\omega}+\eta M\partial_{\rho}}
				\end{aligned}
			\end{align*}	
			is not involutive, we have $n_3=1$. Because of $\C{D_2}=D_1$, item \ref{it:a} is satisfied. In this example, we have
			\begin{align*}
				\begin{aligned}
					D_2^{(1)}&=\Span{\partial_{I_d},\partial_{I_q},\partial_{\rho},\partial_{\psi_d},\partial_{\omega}}=\overline{D}_2\,,
				\end{aligned}
			\end{align*}
			\ie the first derived flag of $D_2$ is already its involutive closure, and thus $n_2=3$. Since $n_2=3$, condition \eqref{eq:compatibilityThm1} of item \ref{it:c} does not occur. Condition \eqref{eq:coupling} is met, in fact we have $\overline{D}_2+[a,D_2]=\mathcal{T}(\mathcal{X})$. Item \ref{it:d} and \ref{it:e} are met with $G_1=\mathcal{T}(\mathcal{X})$. Therefore, according to Theorem \ref{thm:1}, the system is static feedback equivalent to the triangular form \eqref{eq:triangularform}. Since $\Dim{G_1}=\Dim{\overline{D}_2}+1$, the $x_1$-subsystem in a corresponding triangular form \eqref{eq:triangularform} consists only of one integrator chain, since $G_1=\mathcal{T}(\mathcal{X})$, \ie $s=1$, this chain has length one. Thus, according to Section \ref{subse:flatOutputs}, Case 3, flat outputs compatible with the triangular form \eqref{eq:triangularform} are all pairs of functions $(\varphi^1,\varphi^2)$ which satisfy $L^\perp=\Span{\D\varphi^1,\D\Lie_a\varphi^1,\D\varphi^2}=D_2^\perp+\Span{\D\Lie_a\varphi^1}$ and $\Span{\D\varphi^1}=(\overline{D}_2)^\perp$. We have $(\overline{D}_2)^\perp=\Span{\D\theta}$, thus $\varphi^1=\varphi^1(\theta)$. Furthermore, we have $\Lie_a\varphi^1(\theta)=\omega\partial_{\theta}\varphi^1(\theta)$ thus $L^\perp=D_2^\perp+\Span{\D\Lie_a\varphi^1}=\Span{\D\theta,\D\omega,\D\rho}$. Therefore, $\varphi^2=\varphi^2(\theta,\omega,\rho)$, chosen such that $\D\varphi^1\wedge\D\Lie_a\varphi^1\wedge\varphi^2\neq 0$. A possible flat output is thus \eg $\varphi^1=\theta$, $\varphi^2=\rho$. The transformation into the form \eqref{eq:motorTriangular} can be derived systematically by following the six steps of the sufficiency part of the proof of Theorem \ref{thm:1}.
		\vspace{-1.5ex}
	\section{Conclusion}
		\vspace{-2ex}
		We have presented a structurally flat triangular form together with necessary and sufficient conditions for a two-input AI-system to be static feedback equivalent to the proposed triangular form. This provides a sufficient condition for an AI-system to be flat. Future research is devoted to the case where the integrator chains of the $x_3$-subsystem differ in length.
		\vspace{-1ex}
	\bibliography{Bibliography}           

\begin{thebibliography}{27}
\providecommand{\natexlab}[1]{#1}
\providecommand{\url}[1]{\texttt{#1}}
\providecommand{\urlprefix}{URL }
\expandafter\ifx\csname urlstyle\endcsname\relax
  \providecommand{\doi}[1]{doi:\discretionary{}{}{}#1}\else
  \providecommand{\doi}{doi:\discretionary{}{}{}\begingroup
  \urlstyle{rm}\Url}\fi

\bibitem[{Bououden et~al.(2011)Bououden, Boutat, Zheng, Barbot, and
  Kratz}]{BououdenBoutatZhengBarbotKratz:2011}
Bououden, S., Boutat, D., Zheng, G., Barbot, J., and Kratz, F. (2011).
\newblock A triangular canonical form for a class of 0-flat nonlinear systems.
\newblock \emph{International Journal of Control}, 84(2), 261--269.

\bibitem[{Cartan(1914)}]{Cartan:1914}
Cartan, E. (1914).
\newblock Sur l'\'equivalence absolue de certains syst\`emes d'\'equations
  diff\'erentielles et sur certaines familles de courbes.
\newblock \emph{Bulletin de la Soci\'et\'e Math\'ematique de France}, 42,
  12--48.

\bibitem[{Chiasson(1998)}]{Chiasson:1998}
Chiasson, J. (1998).
\newblock A new approach to dynamic feedback linearization control of an
  induction motor.
\newblock \emph{IEEE Trans. Automat. Contr.}, 43, 391--397.

\bibitem[{Fliess et~al.(1992)Fliess, L{\'e}vine, Martin, and
  Rouchon}]{FliessLevineMartinRouchon:1992}
Fliess, M., L{\'e}vine, J., Martin, P., and Rouchon, P. (1992).
\newblock Sur les syst{\`e}mes non lin{\'e}aires diff{\'e}rentiellement plats.
\newblock \emph{Comptes rendus de l'Acad{\'e}mie des sciences. S{\'e}rie I,
  Math{\'e}matique}, 315, 619--624.

\bibitem[{Fliess et~al.(1995)Fliess, L{\'e}vine, Martin, and
  Rouchon}]{FliessLevineMartinRouchon:1995}
Fliess, M., L{\'e}vine, J., Martin, P., and Rouchon, P. (1995).
\newblock Flatness and defect of non-linear systems: introductory theory and
  examples.
\newblock \emph{International Journal of Control}, 61(6), 1327--1361.

\bibitem[{Hunt and Su(1981)}]{HuntSu:1981}
Hunt, L. and Su, R. (1981).
\newblock Linear equivalents of nonlinear time varying systems.
\newblock In \emph{Proceedings 5th International Symposium on Mathematical
  Theory of Networks and Systems (MTNS)}, 119--123.

\bibitem[{Jakubczyk and Respondek(1980)}]{JakubczykRespondek:1980}
Jakubczyk, B. and Respondek, W. (1980).
\newblock On linearization of control systems.
\newblock \emph{Bull. Acad. Polonaise Sci. Ser. Sci. Math.}, 28, 517--522.

\bibitem[{Kolar et~al.(2015)Kolar, Sch{\"o}berl, and
  Schlacher}]{KolarSchoberlSchlacher:2015}
Kolar, B., Sch{\"o}berl, M., and Schlacher, K. (2015).
\newblock Remarks on a triangular form for 1-flat {P}faffian systems with two
  inputs.
\newblock In \emph{Proceedings 1st IFAC Conference on Modelling, Identification
  and Control of Nonlinear Systems (MICNON)}.
\newblock IFAC-PapersOnLine, volume 48, issue 11, pages 109--114.

\bibitem[{Li et~al.(2016)Li, Nicolau, and Respondek}]{LiNicolauRespondek:2016}
Li, S., Nicolau, F., and Respondek, W. (2016).
\newblock Multi-input control-affine systems static feedback equivalent to a
  triangular form and their flatness.
\newblock \emph{International Journal of Control}, 89(1), 1--24.

\bibitem[{Li and Respondek(2012)}]{LiRespondek:2012}
Li, S. and Respondek, W. (2012).
\newblock Flat outputs of two-input driftless control systems.
\newblock \emph{ESAIM: COCV}, 18(3), 774--798.

\bibitem[{Li et~al.(2013)Li, Xu, Su, and Chu}]{LiXuSuChu:2013}
Li, S., Xu, C., Su, H., and Chu, J. (2013).
\newblock Characterization and flatness of the extended chained system.
\newblock In \emph{Proceedings of the 32nd Chinese Control Conference},
  1047--1051.

\bibitem[{Martin and Rouchon(1994)}]{MartinRouchon:1994}
Martin, P. and Rouchon, P. (1994).
\newblock Feedback linearization and driftless systems.
\newblock \emph{Mathematics of Control, Signals and Systems}, 7(3), 235--254.

\bibitem[{Murray(1994)}]{Murray:1994}
Murray, R. (1994).
\newblock Nilpotent bases for a class of non-integrable distributions with
  applications to trajectory generation for nonholonomic systems.
\newblock \emph{Math. Control Signals Systems}, 7, 58--75.

\bibitem[{Nicolau(2014)}]{Nicolau:2014}
Nicolau, F. (2014).
\newblock \emph{Geometry and flatness of control systems of minimal
  differential weight}.
\newblock Theses, INSA de Rouen.

\bibitem[{Nicolau et~al.(2014)Nicolau, Li, and
  Respondek}]{NicolauLiRespondek:2014}
Nicolau, F., Li, S., and Respondek, W. (2014).
\newblock Control-affine systems compatible with the multi-chained form and
  their x-maximal flatness.
\newblock In \emph{Proceedings 21st International Symposium on Mathematical
  Theory of Networks and Systems (MTNS)}, 303--310.

\bibitem[{Nicolau and Respondek(2013)}]{NicolauRespondek:2013}
Nicolau, F. and Respondek, W. (2013).
\newblock Flatness of two-input control-affine systems linearizable via
  one-fold prolongation.
\newblock In \emph{Proceedings 9th IFAC Symposium on Nonlinear Control Systems
  (NOLCOS)}, 493--498.

\bibitem[{Nicolau and Respondek(2014)}]{NicolauRespondek:2014}
Nicolau, F. and Respondek, W. (2014).
\newblock Normal forms for flat control-affine systems linearizable via
  one-fold prolongation.
\newblock In \emph{Proceedings 13th European Control Conference (ECC)},
  2448--2453.

\bibitem[{Nicolau and Respondek(2016{\natexlab{a}})}]{NicolauRespondek:2016-2}
Nicolau, F. and Respondek, W. (2016{\natexlab{a}}).
\newblock Flatness of two-input control-affine systems linearizable via a
  two-fold prolongation.
\newblock In \emph{2016 IEEE 55th Conference on Decision and Control (CDC)},
  3862--3867.

\bibitem[{Nicolau and Respondek(2016{\natexlab{b}})}]{NicolauRespondek:2016}
Nicolau, F. and Respondek, W. (2016{\natexlab{b}}).
\newblock Two-input control-affine systems linearizable via one-fold
  prolongation and their flatness.
\newblock \emph{European Journal of Control}, 28, 20 -- 37.

\bibitem[{Nicolau and Respondek(2017)}]{NicolauRespondek:2017}
Nicolau, F. and Respondek, W. (2017).
\newblock Flatness of multi-input control-affine systems linearizable via
  one-fold prolongation.
\newblock \emph{SIAM J. Control and Optimization}, 55, 3171--3203.

\bibitem[{Nicolau and Respondek(2019)}]{NicolauRespondek:2019}
Nicolau, F. and Respondek, W. (2019).
\newblock Normal forms for multi-input flat systems of minimal differential
  weight.
\newblock \emph{International Journal of Robust and Nonlinear Control}, 29(10),
  3139--3162.

\bibitem[{Nijmeijer and van~der Schaft(1990)}]{NijmeijervanderSchaft:1990}
Nijmeijer, H. and van~der Schaft, A. (1990).
\newblock \emph{Nonlinear Dynamical Control Systems}.
\newblock Springer, New York.

\bibitem[{Schlacher and Sch{\"o}berl(2013)}]{SchlacherSchoberl:2013}
Schlacher, K. and Sch{\"o}berl, M. (2013).
\newblock A jet space approach to check {P}faffian systems for flatness.
\newblock In \emph{Proceedings 52nd IEEE Conference on Decision and Control
  (CDC)}, 2576--2581.

\bibitem[{Sch{\"o}berl(2014)}]{Schoberl:2014}
Sch{\"o}berl, M. (2014).
\newblock \emph{Contributions to the Analysis of Structural Properties of
  Dynamical Systems in Control and Systems Theory - A Geometric Approach}.
\newblock Shaker Verlag, Aachen.

\bibitem[{Sch{\"o}berl et~al.(2010)Sch{\"o}berl, Rieger, and
  Schlacher}]{SchoberlRiegerSchlacher:2010}
Sch{\"o}berl, M., Rieger, K., and Schlacher, K. (2010).
\newblock System parametrization using affine derivative systems.
\newblock In \emph{Proceedings 19th International Symposium on Mathematical
  Theory of Networks and Systems (MTNS)}, 1737--1743.

\bibitem[{Sch{\"o}berl and Schlacher(2014)}]{SchoberlSchlacher:2014}
Sch{\"o}berl, M. and Schlacher, K. (2014).
\newblock On an implicit triangular decomposition of nonlinear control systems
  that are 1-flat - a constructive approach.
\newblock \emph{Automatica}, 50, 1649--1655.

\bibitem[{Silveira et~al.(2015)Silveira, Pereira, and
  Rouchon}]{SilveiraPereiraRouchon:2015}
Silveira, H., Pereira, P., and Rouchon, P. (2015).
\newblock A flat triangular form for nonlinear systems with two inputs:
  Necessary and sufficient conditions.
\newblock \emph{European Journal of Control}, 22, 17 -- 22.

\end{thebibliography}
	\appendix
	\vspace{-1.5ex}
	\section{Supplements}\label{ap:supplements}
		\vspace{-2.1ex}
		In the following, we state a procedure for transforming a driftless system \eqref{eq:driftless} into chained form (or extended chained form if a drift is present), provided it is possible, \ie provided the system meets the corresponding conditions stated in Section \ref{se:knownResults}. According to \cite{LiRespondek:2012} Theorem 2.3, a pair of functions $(\varphi^1,\varphi^2)$ forms a flat output of a driftless system static feedback equivalent to the chained form, if and only if $\D\varphi^1\wedge\D\varphi^2\neq 0$, $L=(\Span{\D\varphi^1,\D\varphi^2})^\perp\subset D^{(n-3)}$ and a regularity condition holds\footnote{The sufficiency part of the proof of this theorem is done constructively and results in a chained form with the components of the flat outputs as top variables. Here, we essentially replicate the sufficiency part of their proof, though argue differently in some places.}. Assume we have a pair of functions $(\varphi^1,\varphi^2)$, which meets the conditions from above (and thus forms a flat output of the considered system). We can transform the system into chained form such that these functions occur as top variables by the following successive procedure. Note that the distribution $D=\Span{b_1,b_2}$ spanned by the input vector fields of a system which is static feedback equivalent to the (extended) chained form, meets the conditions of Lemma \ref{lem:cartan}. We start by straightening out the Cauchy characteristics
		\begin{align}\label{eq:cauchySequence}
			\begin{aligned}
				\C{D^{(1)}}\subset\C{D^{(2)}}\subset\ldots\subset\C{D^{(n-3)}}
			\end{aligned}
		\end{align}
		of the derived flags of $D=\Span{b_1,b_2}$, such that $\C{D^{(i)}}=\Span{\partial_{x^n},\ldots,\partial_{x^{n+1-i}}}$, $i=1,\ldots,n-3$. \mbox{Simultaneously}, we straighten out $L=(\Span{\D\varphi^1,\D\varphi^2})^\perp$ by introducing the components of the flat output as new coordinates, \ie $\bar{x}^j=\varphi^j$, $j=1,2$ (this is always possible since $L\subset D^{(n-3)}$ implies $\C{D^{(n-3)}}\subset L$, see \cite{LiRespondek:2012}). Next, we apply a static feedback to normalize the first equation, \ie $\dot{\bar{x}}^1=\bar{u}^2$. After that, we have a representation of the form
		\begin{align}\label{eq:toChainedForm}
			\hspace{-3ex}\begin{aligned}
				\dot{\bar{x}}^1&=\bar{u}^2\\
				\dot{\bar{x}}^2&=\bar{b}_2^2(\bar{x}^1,\bar{x}^2,x^3)\bar{u}^2\\
				\dot{x}^3&=\bar{b}_2^3(\bar{x}^1,\bar{x}^2,x^3,x^4)\bar{u}^2\\
				&\vdotswithin{=}\\
				\dot{x}^{n-1}&=\bar{b}_2^{n-1}(\bar{x}^1,\bar{x}^2,x^3,\ldots,x^n)\bar{u}^2\\
				\dot{x}^n&=\bar{b}_1^n(\bar{x}^1,\bar{x}^2,x^3,\ldots,x^n)\bar{u}^1+\bar{b}_2^n(\bar{x}^1,\bar{x}^2,x^3,\ldots,x^n)\bar{u}^2.\hspace{-8ex}
			\end{aligned}
		\end{align}
		Next, we successively introduce the functions $\bar{b}_2^i$ as new coordinates from top to bottom, \ie in the first step $\bar{x}^3=\bar{b}_2^2(\bar{x}^1,\bar{x}^2,x^3)$ with all the other coordinates left unchanged, and so on, until in the last step, we complete the transformation to chained form by applying a static feedback to normalize the last equation.
		
		\vspace{-0.5ex}
		The procedure is based on the special structure \eqref{eq:toChainedForm}, which the system takes after applying the above described state and input transformation. In the following, we explain why the system indeed takes this structure. For that, we first show that since Lemma \ref{lem:cartan} applies to the distribution $D$, we have $D^{(i)}=\Span{\bar{b}_2}+\C{D^{(i+1)}}$, $i=0,\ldots,n-4$, \ie the derived flags $D^{(i)}$ are composed of the one-dimensional distribution spanned by the vector field $\bar{b}_2$ and the Cauchy characteristic distributions of their next derived flags (from which it also follows that the input $\bar{u}^1$ can only occur in the very last equation of \eqref{eq:toChainedForm}). To show this, note that $\bar{b}_2$ has a component in the $\partial_{\bar{x}^1}$-direction. Thus, it cannot belong to any of the Cauchy characteristics $\C{D^{(i+1)}}=\Span{\partial_{x^n},\ldots,\partial_{x^{n-i}}}$, $i=0,\ldots,n-4$. However, because of $\bar{b}_2\in D$, it belongs to all the derived flags $D^{(i)}$. Furthermore, because of Lemma \ref{lem:cartan}, we have $\Dim{\C{D^{(i+1)}}}=\Dim{D^{(i)}}-1$, $i=0,\ldots,n-4$. Thus, $\Span{\bar{b}_2}$ completes $\C{D^{(i+1)}}$ to $D^{(i)}$. Furthermore, by construction we have $L=(\Span{\D\varphi^1,\D\varphi^2})^\perp\subset D^{(n-3)}$ and we also have $\bar{b}_2\notin L$. Thus, $D^{(n-3)}=\Span{\bar{b}_2}+L$. In conclusion, we therefore have $D^{(i)}=\Span{\partial_{x^n},\ldots,\partial_{x^{n-i}},\bar{b}_2}$, $i=0,\ldots,n-3$ and thus, $\partial_{x^j}\bar{b}_2^k=0$ for $k+2\leq j\leq n$, $k=2,\ldots,n-2$, and $\partial_{x^j}\bar{b}_2^k\neq 0$ for $j=k+1$, $k=2,\ldots,n-1$ follows. This exactly describes the triangular dependence of the functions $\bar{b}_2^i$, $i=2,\ldots,n-1$ on the states of the system, which enables us to successively introduce these functions as new states.
\end{document}